\documentclass[reqno]{amsart}
\usepackage{hyperref}
\usepackage{mathrsfs}

\AtBeginDocument{}

\begin{document}
\title[\hfil Regularity  of pullback attractors]
{  Regularity  of pullback attractors for non-autonomous stochastic
 FitzHugh-Nagumo systems with additive noises on unbounded domains}

\author{Wenqiang Zhao}  

\address{Wenqiang Zhao \newline
School Of Mathematics and Statistics, Chongqing Technology and Business
University,  Chongqing 400067, China}
\email{gshzhao@sina.com}

\subjclass[2000]{60H15, 35R60, 35B40, 35B41}
\keywords{ Random dynamical systems;
 non-autonomous FitzHugh-Nagumo systems; upper semi-continuity; $\mathcal{D}$-pullback attractor; multiplicative noise}

\begin{abstract}
In this paper, we prove the existences of pullback attractors  in $L^{p}(\mathbb{R}^N)\times L^{2}(\mathbb{R}^N)$
for stochastic Fitzhugh-Nagumo system driven by both additive noises and  deterministic non-autonomous forcings. The nonlinearity is polynomial like growth with  exponent $p-1$.  The asymptotic compactness for the cocycle in $L^{p}(\mathbb{R}^N)\times L^{2}(\mathbb{R}^N)$ is  proved  by using asymptotic a priori method, where the plus and minus signs of the nonlinearity at large value are not required.
\end{abstract}

\maketitle
\numberwithin{equation}{section}
\newtheorem{theorem}{Theorem}[section]
\newtheorem{lemma}[theorem]{Lemma}
\newtheorem{remark}[theorem]{Remark}
\newtheorem{definition}[theorem]{Definition}
\allowdisplaybreaks

\section{Introduction}

In this paper, we consider the regularity  of pullback attractors for
the following non-autonomous FitzHugh-Nagumo system defined on $\mathbb{R}^N$ perturbed by additive noises:
\begin{equation}
\begin{cases}
d\tilde{u}+(\lambda\tilde{ u}-\Delta\tilde{ u}+\alpha \tilde{v})dt=f(x,\tilde{u})dt+g(t,x)dt+ h_1 d\omega_1(t),\\
d\tilde{v}+(\sigma \tilde{v}-\beta\tilde{u})dt=h(t,x)dt+h_2d\omega_2(t),\\
\tilde{u}(x,\tau)=\tilde{u}_0(x),\\
\tilde{v}(x,\tau)=\tilde{v}_0(x),
\end{cases}
\end{equation}
where the initial condition $(\tilde{u}_0,\tilde{v}_0)\in L^2(\mathbb{R}^N)\times L^2(\mathbb{R}^N)$, the coefficients $\lambda,\alpha,\beta,\sigma$ are  positive constants, $h_1$
and $h_2$ are given functions on $\mathbb{R}^N$ satisfying some regular conditions,
the non-autonomous terms $g,h\in L^2_{loc}(\mathbb{R}, L^2(\mathbb{R}^N))$, the nonlinearity $f$ is a smooth function satisfying some polynomial growth with exponent $p-1,p>2$ and $\omega(t)=(\omega_1(t),\omega_2(t))$ is a Wiener process defined on a probability space $(\Omega,
\mathcal {F}, {P})$, which will be specified later.

The deterministic FitzHugh-Nagumo system is an important mathematical model to describe the signal transmission across axons in neurobiology, which is well studied  in the literature, see, e.g., \cite{Nagumo,FitzHugh,Jones, Sys1,Sys2,Sys3}.
In the stochastic case,
 Wang \cite{Wang2} proved the existence and uniqueness of random attractors in the initial space $L^2(\mathbb{R}^N)\times L^2(\mathbb{R}^N)$ for the autonomous FitzHugh-Nagumo system. For the general non-autonomous forcings $g$ and $h$,  Adili and Wang \cite{Adili1, Adili2} obtained  the pullback attractors in the initial space $L^2(\mathbb{R}^N)\times L^2(\mathbb{R}^N)$, and Bao \cite{Bao} strengthened  this result and proved the regularity of pullback attractors in the non-initial space $H^1(\mathbb{R}^N)\times L^2(\mathbb{R}^N)$.  Recently, \cite{Huang,Gu,Gu1} studied the existences of random attractors for stochastic lattice FitzHugh-Nagumo system

However, to our knowledge, there are no literature to investigate the high-order integrability of solutions  of problem (1.1), even for the deterministic case.

In this paper,  we prove that  the problem (1.1)  admits  a unique pullback
attractor in the non-initial space $L^\varpi(\mathbb{R}^N)\times L^2(\mathbb{R}^N)$ for every $\varpi\in(2,p]$, with the functions $f,g $ and $h$  satisfying almost the same conditions  as in \cite{Adili2}. Actually, we  show that the obtained  pullback attractors in a family of spaces $L^\varpi(\mathbb{R}^N)\times L^2(\mathbb{R}^N)(\varpi\in(2,p])$  are the same object as
 in the initial space $L^2(\mathbb{R}^N)\times L^2(\mathbb{R}^N)$. Furthermore,  it
is compact and attracts every nonempty bounded subset of $L^2(\mathbb{R}^N)\times L^2(\mathbb{R}^N)$
under the topology of $L^\varpi(\mathbb{R}^N)\times L^2(\mathbb{R}^N)(\varpi\in(2,p])$.  The asymptotic compactness in $L^p(\mathbb{R}^N)\times L^2(\mathbb{R}^N)$ for the cocycle  is
 proved by using asymptotic a priori method which is developed by \cite{Zhong}.  Because of the interactions of the variables in  the system (1.1),
 proving that the unbounded part of solutions vanishes in $L^p$ norm is a more technical task.  We prove this in a relatively simple way just using asymptotic a priori method, in which  the plus and minus signs of the nonlinearity at large value are not required, comparing with some similar work as in \cite{Liyangrong1,Zhao4,Lijia,Yin}.

This paper is organized as follows. In the next section, we recall some preliminaries required for our further discussions and then establish an abstract result regarding
the existence of pullback attractors in some non-initial spaces,
which can be applicable to stochastic partial differential equations (SPDEs) with both random noise and a deterministic non-autonomous term. In Section 3, we give the assumptions on $g, h$ and $f$ and define a continuous cocycle for problem (1.1).
Finally we prove the existences of pullback attractors in $L^\varpi(\mathbb{R}^N)\times L^2(\mathbb{R}^N)$ in Section 4.

\section{Preliminaries and abstract results}

Let $(X, \|.\|_X)$ and $(Y, \|.\|_Y)$  be two complete separable Banach spaces  with Borel
sigma-algebra $\mathcal {B}(X)$ and $\mathcal {B}(Y)$, respectively. For convenience, we call $X$ \emph{the initial space} (contains all initial data of a SPDE)
  and $Y$ \emph{the associated non-initial space} (usually the  regular solutions (of a SPDE) located space).

In this section,  we give a sufficient standard for the existence of pullback attractors
in the non-initial  space $Y$ for
the random dynamical system (RDS) over two parametric spaces.  For the existence of random attractors in the non-initial  space $Y$ for the RDS over one parametric space,  the readers may refer to \cite{Zhao0,Zhao3,Zhao4,Lijia,Liyangrong1,Liyangrong2,Liyangrong3,Yin}.

We also mention that regarding the existence of random attractors  in the initial space $X$ for the RDS over one parametric space, the good references
 are \cite{Arn,Bates,Rand1,Rand2,Rand3,Rand4}. However, here we recall from  \cite{Wang4} some basic notions regarding RDS over two parametric spaces.
Let $\Omega_1$ be a nonempty set on which there is a mapping $\vartheta_{1,t}: \Omega_1\rightarrow \Omega_1$  such that $\vartheta_{1,0}$ is the identity on $\Omega_1$
and $\vartheta_{1,s+t}=\vartheta_{1,t}\circ \vartheta_{1,s}$ for all $s,t\in\mathbb{R}$. Let $(\Omega_2,\mathcal{F}_2,{P})$  be  a probability space on which there is  a measure preserving transformation $\vartheta_{2,t}$ such that $\vartheta_{2,0}$ is the identity on $\Omega_2$
and $\vartheta_{2,s+t}=\vartheta_{2,t}\circ \vartheta_{2,s}$ for all $s,t\in\mathbb{R}$. We often use $(\Omega_1,\{\vartheta_{1,t}\}_{t\in\mathbb{R}})$ and $(\Omega_2,\mathcal{F}_2,{P},\{\vartheta_{2,t}\}_{t\in\mathbb{R}})$ to denote the two parametric spaces respectively.

A mapping
 $\varphi: \mathbb{R}^+\times\Omega_1\times \Omega_2\times X\mapsto X$ is called  a cocycle on $X$ over $(\Omega_1,\{\vartheta_{1,t}\}_{t\in\mathbb{R}})$ and $(\Omega_2,\mathcal{F}_2,{P},\{\vartheta_{2,t}\}_{t\in\mathbb{R}})$
if for every $\omega_1\in \Omega_1, \omega_2\in\Omega_2$ and $t,s\in \mathbb{R}^+$, the following conditions are satisfied:

(i) $\varphi(.,\omega_1,.,.):\mathbb{R}^+\times \Omega_2\times X\mapsto X$ is $(\mathcal{B}(\mathbb{R}^+)\times \mathcal{F}_2\times\mathcal{B}(X),\mathcal{B}(X))$-measurable,

(ii) $\varphi(0, \omega_1,\omega_2,.)$ is the identity on X,

(iii) $\varphi(t+s,\omega_1, \omega_2,.)=\varphi(t,\vartheta_{1,s}\omega_1, \vartheta_{2,s}\omega_2, .)\circ \varphi(s,\omega_1, \omega_2, .).$\\
In addition, if $\varphi(t,\omega_1, \omega_2,.): X\rightarrow X$ is continuous for every $t\in
\mathbb{R}^+, \omega_1\in \Omega_1, \omega_2\in\Omega_2$, then $\varphi$ is called a continuous cocycle on $X$ over $(\Omega_1,\{\vartheta_{1,t}\}_{t\in\mathbb{R}})$ and $(\Omega_2,\mathcal{F}_2,{P},\{\vartheta_{2,t}\}_{t\in\mathbb{R}})$.

A set-valued mapping $K: \Omega_1\times \Omega_2\mapsto 2^X$ is called measurable in $X$ with respect to $\mathcal{F}_2$ in $\Omega_2$ if the mapping $\omega_2\in\Omega_2\mapsto \mbox{dist}_X(x,K(\omega_1,\omega_2))$ is ($\mathcal{F}_2,\mathcal{B}(\mathbb{R}))$-measurable for every fixed $x\in X$
and $\omega_1\in\Omega_1$, where $\mbox{dist}_X$ is the Haustorff semi-metric in $X$. And furthermore if the value $K(\omega_1,\omega_2)$ is a closed nonempty subset of $X$ for every $\omega_1\in\Omega_1$ and $\omega_2\in\Omega_2$ then $\{K(\omega_1,\omega_2);\omega_1\in\Omega_1,\omega_2\in\Omega_2\}$ is called a closed measurable subset of $X$ with respect to $\mathcal{F}_2$ in $\Omega_2$.

In the sequel, we assume that the cocycle $\varphi$ together with $X$ and $Y$ satisfy the following\\

\emph{\textbf{(H1)£º}\ \  For every fixed $t>0, \omega_1\in\Omega_1$ and $\omega_2\in\Omega_2$, $\varphi(t,\omega_1,\omega_2, .):X\mapsto   Y$. }\\

We always assume that $\varphi$ is a continuous cocycle on the initial space $X$ over $(\Omega_1,\{\vartheta_{1,t}\}_{t\in\mathbb{R}})$ and $(\Omega_2,\mathcal{F}_2,{P},\{\vartheta_{2,t}\}_{t\in\mathbb{R}})$,
and $\mathcal{D}$ a collection of some families of nonempty subsets of $X$ parametrized by $\omega_1\in\Omega_1$ and $\omega_2\in\Omega_2$, \emph{i.e.,}
$$
\mathcal{D}=\{D=\{D(\omega_1,\omega_2)\in 2^X;D(\omega_1,\omega_2)\neq\emptyset,\omega_1\in\Omega_1, \omega_2\in\Omega_2 \};f_D\ \mbox{satisfies some conditions}\}.
$$
Let $D_1$ and $D_2$ belong to $\mathcal{D}$. We say $ D_1=D_2$ if and only if $D_1(\omega_1, \omega_2)=D_2(\omega_1, \omega_2) $ for all $\omega_1\in\Omega_1$ and $\omega_2\in\Omega_2$.
Note that for our problem the  continuity of the cocycle $\varphi$ in the non-initial space $Y$ is not necessary.
  \\

\textbf{Definition 2.1.} \emph{Let  $\mathcal{D}$ be  a  collection  of some families of nonempty subsets of $X$. $ \{K(\omega_1,\omega_2);\omega_1\in\Omega_1,\omega_2\in\Omega_2\}\in\mathcal{D}$
 is called a $\mathcal{D}$-pullback absorbing set for
a cocycle  $\varphi$ in $X$ if for every $\omega_1\in\Omega_1, \omega_2\in\Omega_2$ and  $D\in \mathcal{D}$,
 there exists a absorbing time $T=T(\omega_1,\omega_2,D)>0$ such that
$$\varphi (t,\vartheta_{1,-t}\omega_1, \vartheta_{2,-t}\omega_2, D(\vartheta_{1,-t}\omega_1,\vartheta_{2,-t}\omega_2))\subseteq K(\omega_2, \omega_2)\ \ \ \ \ \ \mbox{for all}\ t\geq
T.$$}

 \textbf{Definition 2.2.} \emph{Let  $\mathcal{D}$ be  a  collection  of some families of nonempty subsets of $X$. A cocycle  $\varphi$ is said to be
$\mathcal{D}$-pullback asymptotically compact in $X\ (\mbox{resp. in}\ Y)$ if for every $\omega_1\in\Omega_1$ and $\omega_2\in\Omega_2$,
$$\{\varphi(t_n, \vartheta_{1,-t_n}\omega_1, \vartheta_{2,-t_n}\omega_2, x_n)\}\ \mbox{has a
convergent subsequence in}\ X(\mbox{resp. in}\ Y)$$
whenever $t_n\rightarrow \infty$ and
$x_n\in D(\vartheta_{1,-t}\omega_1, \vartheta_{2,-t}\omega_2)$ with $D=\{D(\omega_1, \omega_2);\omega_1\in \Omega_1,\omega_2\in \Omega_2\} \in \mathcal{D}$.}\\

\textbf{Definition 2.3.} \emph{Let  $\mathcal{D}$ be  a  collection  of some families of nonempty subsets of $X$. $ \{\mathcal{A}(\omega_1,\omega_2);\omega_1\in\Omega_1,\omega_2\in\Omega_2\}\in\mathcal{D}$ is called  a $\mathcal{D}
$-pullback attractor for a cocycle $\varphi$ in $X\ (\mbox{resp. in}\ Y)$ over $(\Omega_1,\{\vartheta_{1,t}\}_{t\in\mathbb{R}})$ and $(\Omega_2,\mathcal{F}_2,{P},\{\vartheta_{2,t}\}_{t\in\mathbb{R}})$ if the following statements are satisfied:}

\emph{(i)  $\mathcal {A}$ is measurable in $X$ with respect to the $P$-completion  of  $\mathcal{F}_2$ in $\Omega_2$  and $\mathcal {A}(\omega_1,\omega_2)$ is compact in $X\ (\mbox{resp. in}\ Y)$ for all $\omega_1\in \Omega_1, \omega_2\in\Omega_2$,}

 \emph{(ii) $\mathcal {A}$ is invariant, that is, for every $\omega_1\in \Omega_1, \omega_2\in\Omega_2$,
 $\varphi(t,\omega_1,\omega_2,\mathcal{A}(\omega_1,\omega_2))=\mathcal{A}(\vartheta_{1,t}\omega_1,\vartheta_{2,t}\omega_2)$ for all $t\geq0,$}

\emph{(iii) $\mathcal {A}$  attracts every element $D=\{D(\omega_1, \omega_2);\omega_1\in \Omega_1,\omega_2\in \Omega_2\} \in \mathcal{D}$ in $X\ (\mbox{resp. in}\ Y)$, that is, for every $\omega_1\in \Omega_1, \omega_2\in\Omega_2$,
$$\mathop {\lim }\limits_{t\rightarrow +\infty}\mbox{dist}_X(\varphi(t,\vartheta_{1,-t}\omega_1,\vartheta_{2,-t}\omega_2, D(\vartheta_{1,-t}\omega_1, \vartheta_{2,-t}\omega_2)),\mathcal{A}(\omega_1,\omega_2))= 0$$
$$(\mbox{resp.}\mathop {\lim }\limits_{t\rightarrow +\infty}\mbox{dist}_Y(\varphi(t,\vartheta_{1,-t}\omega_1,\vartheta_{2,-t}\omega_2, D(\vartheta_{1,-t}\omega_1, \vartheta_{2,-t}\omega_2)),\mathcal{A}(\omega_1,\omega_2))= 0).$$
}

 The following Theorem 2.6 is a useful tool to prove the existence of pullback attractors for SPDEs with both non-autonomous deterministic term
 and random term in a associated non-initial space. We need to assume that the spaces
 $X$ and $Y$ satisfy sequence limits uniqueness property:\\

\emph{\textbf{(H2)}\ \  If $\{x_n\}_n\subset X\cap Y$ such that $x_n\rightarrow x$ in
$X$ and $x_n\rightarrow y$ in $Y$, respectively, then we have $x=y$.}\\

 \textbf{ Theorem  2.6.} \emph{Let  $\mathcal{D}$ be  a  collection  of some families of nonempty subsets of $X$ which is inclusion closed. Suppose that $\varphi$ is a continuous cocycle on the initial space $X$ over $(\Omega_1,\{\vartheta_{1,t}\}_{t\in\mathbb{R}})$ and $(\Omega_2,\mathcal{F}_2,{P},\{\vartheta_{2,t}\}_{t\in\mathbb{R}})$. Assume that}

\emph{(i) $\varphi$ has  a closed and measurable (w.r.t. the $P$-completion of $\mathcal{F}_2$ in $\Omega_2$) $\mathcal{D}$-pullback absorbing set
$K=\{K(\omega_1,\omega_2);\omega_1\in\Omega_1,\omega_2\in\Omega_2\}\in \mathcal{D}$ in $X$;}

\emph{(ii) $\varphi$ is  $\mathcal{D}$-pullback asymptotically compac}t in $X$.
\emph{Then the cocycle  $\varphi$ has a unique $\mathcal{D}$-pullback  attractor $\mathcal{A}_X=\{\mathcal{A}_X(\omega_1,\omega_2);\omega_1\in\Omega_1,\omega_2\in\Omega_2\}$ in $X$,  which is given by
\begin{align}\label{ff01}
\mathcal{A}_X(\omega_1,\omega_2)=\cap_{s>0}\overline{\cup_{t\geq s} \varphi(t,\vartheta_{1,-t}\omega_1,\vartheta_{2,-t}\omega_2, K(\vartheta_{1,-t}\omega_1,\vartheta_{2,-t}\omega_2))}^X, \ \
\omega_1\in\Omega_1,\omega_2\in\Omega_2,
 \end{align} where the closure is taken over in $X$.}

\emph{If further
(H1)-(H2) hold and  $\varphi$ is  $\mathcal{D}$-pullback asymptotically compact in $Y$,}
\emph{ then the cocycle  $\varphi$ has a unique $\mathcal{D}$-pullback  attractor $\mathcal{A}_Y=\{\mathcal{A}_Y(\omega_1,\omega_2);\omega_1\in\Omega_1,\omega_2\in\Omega_2\}\in \mathcal{D}$ in $Y$,
which is given by
\begin{align}\label{ff02}
\mathcal{A}_Y(\omega_1,\omega_2)=\cap_{s>0}\overline{\cup_{t\geq s} \varphi(t,\vartheta_{1,-t}\omega_1,\vartheta_{2,-t}\omega_2, K(\vartheta_{1,-t}\omega_1,\vartheta_{2,-t}\omega_2))}^Y, \ \
\omega_1\in\Omega_1,\omega_2\in\Omega_2,
 \end{align}where the closure is taken over in $Y$.}
\emph{In addition, we have $\mathcal{A}_Y=\mathcal{A}_X$, i.e., for every $\omega_1\in\Omega_1,\omega_2\in\Omega_2$, $\mathcal{A}_Y(\omega_1,\omega_2)=\mathcal{A}_X(\omega_1,\omega_2)$}.\\

{\bf Proof} \ \  The first result is stated in \cite{Wang4}. We only prove the second result.
Indeed, (\ref{ff02}) makes sense by using \emph{(H1)}, and$\mathcal{A}_Y\neq\emptyset$ by the asymptotic compactness of the cocycle $\varphi$ in $Y$. In the following, we  show that  $\mathcal{A}_Y$
is a unique $\mathcal{D}$-pullback attractor in $Y$.

 \emph{Step 1.} Since $\mathcal{A}_X$ is measurable in $X$ (with respect to the $P$-completion of $\mathcal{F}_2$ in $\Omega_2$)  and $\mathcal{A}_X\in\mathcal{D}$ is invariant (the measurability of $\mathcal{A}_X$ is proved by Theorem 2.14 in \cite{Wang4}), then
 the measurability and invariance of  $\mathcal{A}_Y$  can be proved by showing that $\mathcal{A}_Y=\mathcal{A}_X$.   We have to prove that for every $\omega_1\in\Omega_1,\omega_2\in\Omega_2$, $\mathcal{A}_Y(\omega_1,\omega_2)=\mathcal{A}_X(\omega_1,\omega_2)$, where $\mathcal{A}_X(\omega_1,\omega_2)$ and $\mathcal{A}_Y(\omega_1,\omega_2)$ defined in (\ref{ff01}) and (\ref{ff02}) respectively. Indeed,  taking $x\in \mathcal{A}_X(\omega_1,\omega_2)$, by (\ref{ff01}), there are $t_n\rightarrow+\infty$
and $x_n\in K(\vartheta_{1,-t_n}\omega_1,\vartheta_{2,-t_n}\omega_2)$ such that
\begin{align} \label{ff3}
 \varphi(t_n, \vartheta_{1,-t_n}\omega_1,\vartheta_{2,-t_n}\omega_2, x_n)\xrightarrow[n\rightarrow\infty]{\|.\|_X}x.
 \end{align}
Since $\varphi$ is
$\mathcal{D}$-asymptotically compact in $Y$, then there is a $y\in Y$ such that up to a subsequence (we still denote this by itself),
\begin{align}\label{ff4}
 \varphi(t_n, \vartheta_{1,-t_n}\omega_1,\vartheta_{2,-t_n}\omega_2, x_n)\xrightarrow[n\rightarrow\infty]{\|.\|_Y}y.
 \end{align} Then it follows from (\ref{ff02}) that $y\in \mathcal{A}_Y(\omega_1,\omega_2)$. Then by \emph{(H2)} and (\ref{ff3}) and (\ref{ff4}) we have $x=y\in \mathcal{A}_X(\omega_1,\omega_2)$ and thus
 $ \mathcal{A}_X(\omega_1,\omega_2)\subseteq  \mathcal{A}_Y(\omega_1,\omega_2)$. The inverse inclusion can be proved similarly. Thus $\mathcal{A}_X=\mathcal{A}_Y$ as required.

 \emph{Step 2.}  We prove the attraction of  $\mathcal{A}_Y$
in $Y$ by a contradiction argument. Indeed,  if for fixed $\omega_1\in\Omega_1$ and $\omega_2\in\Omega_2$,
there exist $\delta>0$, $x_n\in D(\vartheta_{1,-t_n}\omega_1,\vartheta_{2,-t_n}\omega_2)$ with $D=\{D(\omega_1, \omega_2);\omega_1\in \Omega_1,\omega_2\in \Omega_2\} \in\mathcal{D}$ and
$t_n\rightarrow +\infty$  such that
\begin{align}\label{ff5}
 \mbox{dist}_Y\Big(\varphi(t_n,\vartheta_{1,-t_n}\omega_1,\vartheta_{2,-t_n}\omega_2,x_n),
 \mathcal{A}_Y(\tau,\omega)\Big)\geq \delta.
 \end{align}
By the asymptotic compactness of $\varphi$ in $Y$, there exists $y_0\in Y$  such that up to a
subsequence,
\begin{align}\label{ff6}
 \varphi(t_n,\vartheta_{1,-t_n}\omega_1,\vartheta_{2,-t_n}\omega_2,x_n)\xrightarrow[n\rightarrow\infty]{\|.\|_Y}y_0.
 \end{align}
On the other hand, by our assumption (i), there exist a large time $T>0$ such that
\begin{align}\label{ff7}
 y_n=\varphi(T,\vartheta_{1,-t_n}\omega_1,\vartheta_{2,-t_n}\omega_2,x_n)&=\varphi(T,\vartheta_{1,-T}\vartheta_{1,-t_n+T}\omega_1,\vartheta_{2,-T}\vartheta_{2,-t_n+T}\omega_2,x_n)\notag\\
 &\in K(\vartheta_{1,-t_n+T}\omega_1,\vartheta_{2,-t_n+T}\omega_2).
 \end{align}
Then by the cocycle property of $\varphi$, along with (\ref{ff6}) and (\ref{ff7}), we infer that as $n\rightarrow\infty$,
\begin{align}
 \varphi(t_n,\vartheta_{1,-t_n}\omega_1,\vartheta_{2,-t_n}\omega_2,x_n)&=\varphi(t_n-T,\vartheta_{1,-t_n+T}\omega_1,\vartheta_{2,-t_n+T}\omega_2,y_n)\rightarrow y_0\ \  \ \mbox{in}\ Y.\notag
 \end{align} Therefore by (\ref{ff02}), $y_0\in \mathcal{A}_Y(\omega_1,\omega_2)$. This implies that
 \begin{align}
 \mbox{dist}_Y\Big(\varphi(t_n,\vartheta_{1,-t_n}\omega_1,\vartheta_{2,-t_n}\omega_2,x_n),
 \mathcal{A}_Y(\omega_1,\omega_2)\Big)\rightarrow 0\notag
 \end{align} as $n\rightarrow\infty$, which is a contradiction to (\ref{ff5}).

\emph{Step 3.} It remains to prove the compactness of
$A_Y$ in $Y$.  Let
$\{y_n\}_{n=1}^\infty$  be a sequence in
$A_Y(\omega_1,\omega_2)$. By the invariance of $A_Y(\omega_1,\omega_2)$ which is proved in Step 1, we
 have
 $$
 \varphi(t,\vartheta_{1,-t}\omega_1,\vartheta_{2,-t}\omega_2, \mathcal{A}_Y(\vartheta_{1,-t}\omega_1,\vartheta_{2,-t}\omega_2))=\mathcal{A}(\omega_1,\omega_2).
 $$
 Then it follows that there is a sequence $\{z_n\}_{n=1}^\infty\in \mathcal{A}_Y(\vartheta_{1,-t}\omega_1,\vartheta_{2,-t}\omega_2)$
such that for every $n\in\mathbb{N}$,
 \begin{align}
 y_n=\varphi(t_n,\vartheta_{1,-t_n}\omega_1,\vartheta_{2,-t_n}\omega_2,z_n).\notag
 \end{align}
Note that $A_Y\in \mathcal{D}$. Then by the asymptotic compactness of $\varphi$ in $Y$, $\{y_n\}$  has a convergence subsequence in $Y$, \emph{i.e.},
 there is a $y_0\in Y$ such that
 $$
 \lim_{n\rightarrow\infty}y_n=y_0\ \ \mbox{in}\ Y.
 $$ But  $A_Y(\omega_1,\omega_2)$ is closed in $Y$, so $y_0\in A_Y(\omega_1,\omega_2)$.

 The uniqueness is easily followed by the attraction property of $\varphi$ and $A_Y\in \mathcal{D}$. This completes the total proofs.$\ \ \ \ \ \
\Box$\\

\textbf{Remark 2.7.} (i) We emphasize that the assumption \emph{(H1)} is necessary to guarantee that the closure of the set $\varphi(t,\vartheta_{1,-t}\omega_1,\vartheta_{2,-t}\omega_2,K(\vartheta_{1,-t}\omega_1,\vartheta_{2,-t}\omega_2))$ in $Y$ makes sense for all $t>0$, as in (\ref{ff02}).

(ii) We emphasize that the attractor $A_Y$ is measurable in $X$ and completely structured by a closed and measurable  $\mathcal{D}$-pullback absorbing set in the initial space $X$. Thus the existence of absorbing set in the non-initial space $Y$ is unnecessary. This is different from the construction in \cite{Bao}.

(iii ) The above Theorem 2.6 generalizes  Theorem 3.1  in \cite{Liyangrong3} regarding
the bi-spatial random attractor which is suitable for the random dynamical system over one parametric space. Our result here is applicable to
the SPDE with random noises and deterministic non-autonomous forcing.\\

In particular, if the initial space $X=L^2(\mathbb{R}^N)$ and the associated non-initial space $Y=L^r(\mathbb{R}^N)$, then we have an easy condition on the
asymptotic compactness  of the cocycle $\varphi$ in the non-initial space $Y=L^r(\mathbb{R}^N),r>2$. First we let $\mathcal{D}$ be the collection of some family of nonempty
subsets of $L^2(\mathbb{R}^N)$. Then we have\\

\textbf{Proposition 2.8.} \emph{Assume that  the cocycle $\varphi$ is asymptotically compact in $L^2(\mathbb{R}^N)$. Given $\omega_1\in\Omega_1, \omega_2\in\Omega_2$ and $D=\{D(\omega_1, \omega_2);\omega_1\in \Omega_1,\omega_2\in \Omega_2\}\in\mathcal{D}$, suppose that for every $\eta>0$, there are constants $M=M(\omega_1,\omega_2,D,\eta)>0$ and $T=T(\omega_1,\omega_2,D)>0$ such that
\begin{align}\label {}
\sup\limits_{t\geq T}\sup\limits_{u_0\in{D}(\vartheta_{1,-t}\omega_1,\vartheta_{2,-t}\omega_2))}\int\limits_{\mathbb{R}^N(|\varphi_t|\geq M)}|\varphi_t|^rdx\leq\eta,\notag
\end{align} where $\varphi_t=\varphi(t,\vartheta_{1,-t}\omega_1,\vartheta_{2,-t}\omega_2, u_0)$.
Then the cocycle $\varphi$ is asymptotically compact in $L^r(\mathbb{R}^N)$, that is, for every $\omega_1\in\Omega_1$ and $\omega_2\in\Omega_2$, the sequence $\{\varphi(t, \vartheta_{1,-t_n}\omega_1,\vartheta_{2,-t_n}\omega_2, {u}_{0,n})\}$ has a convergent subsequence in
 $L^r(\mathbb{R}^N)$ whenever $t_n\rightarrow+\infty$ and ${u}_{0,n}\in D(\vartheta_{1,-t_n}\omega_1,\vartheta_{2,-t_n}\omega_2)\in \mathcal{D}$.}\\

\emph{Proof}\ \  The proof is very similar to the one parameter case, in which we can find the proof in \cite{Zhao3,Liyangrong1,Liyangrong2} and so we omit it here.
 $\ \ \  \ \ \ \ \ \ \Box$\\

\section{ Non-autonomous FitzHugh-Nagumo system on $\mathbb{R}^N$ with additive noises}

In this section, we give  the conditions on the functions $f,g $ and $h$. We also give the parametric dynamical systems $(\Omega_1,\{\vartheta_{1,t}\}_{t\in\mathbb{R}}$ and $(\Omega_2,\mathcal{F}_2,P,\{\vartheta_{2,t}\}_{t\in\mathbb{R}})$ respectively,
and show the existence of continuous cocycle  in the initial space $L^2(\mathbb{R}^N)\times L^2(\mathbb{R}^N)$.
For the non-autonomous FitzHugh-Nagumo system (1.1), the nonlinearity $f(x,s)$
satisfies almost the same assumptions as in \cite{Adili2}, \emph{i.e.},
for all $ x\in\mathbb{R}^N$ and $s\in\mathbb{R}$,
\begin{align}\label {a1}
&f(x,s)s\leq -\alpha_1|s|^p+\psi_1(x),
\end{align}
\begin{align}\label {a2}
&|f(x,s)|\leq \alpha_2|s|^{p-1}+\psi_2(x),
\end{align}
\begin{align}\label {a3}
&\frac{\partial f}{\partial s}f(x,s)\leq \alpha_3 ,
\end{align}
\begin{align}\label {a4}
&\Big|\frac{\partial f}{\partial x}f(x,s)\Big|\leq \psi_3(x).
\end{align} where $p>2$, $\alpha_i> 0(i=1,2,3)$ are determined constants,  $\psi_1\in L^{1}(\mathbb{R}^N)\cap L^{\frac{p}{2}}(\mathbb{R}^N)$,
$\psi_2$ and $\psi_3\in L^2(\mathbb{R}^N)$. The non-autonomous terms
 $g\in L^2_{loc}(\mathbb{R}, L^2(\mathbb{R}^N))$ and $h\in L^2_{loc}(\mathbb{R}, L^2(\mathbb{R}^N))$  satisfy that for every $\tau\in\mathbb{R}$,
\begin{align}\label {a5}
\int\limits_{-\infty}^\tau e^{\delta s} (\|g(s,.)\|_{L^2(\mathbb{R}^N)}^2+\|h(s,.)\|_{L^2(\mathbb{R}^N)}^2)ds<+\infty.
\end{align} where $\delta=\min\{\lambda,\sigma\}$, $\lambda $ and $\sigma$ are as in (1.1). The functions $h_1\in H^2(\mathbb{R}^N)\cap W^{2,p}(\mathbb{R}^N)\cap L^{\infty}(\mathbb{R}^N)$
and $h_2\in H^1(\mathbb{R}^N)$. We remark that the assumption $\psi_1$ is essentially bounded is not required and then the plus and minus signs of the
nonlinearity $f(x,s)$ at large value are unknown by condition (\ref{a1}).

For the probability space $(\Omega,\mathcal{F},P)$, we write
 $\Omega=\{\omega\in
C(\mathbb{R},\mathbb{R}^2); \omega(0) =0\}$. Let
 $\mathcal {F}$ be the
 Borel $\sigma$-algebra induced by the compact-open topology of
$\Omega$ and ${P}$ be  the corresponding
 Wiener measure on $(\Omega,\mathcal{F})$.
  We define a shift operator $\vartheta$ on $\Omega$ by
$$
\vartheta_t\omega(s)=\omega(s+t)-\omega(t), \ \mbox{for every}\ \omega\in\Omega, t,s\in\mathbb{R}.
$$
It is easy to show that
 $\vartheta_t$ is a measure preserving transformation on the given probability space $(\Omega,\mathcal{F},P)$.
Let  $\{\vartheta_{1,t}\}_{t\in\mathbb{R}}$ be a group on $\mathbb{R}$ given by  $\vartheta_{1,t}(\tau)=\tau+t$ for all $\tau,t\in \mathbb{R}$.  Then $(\Omega_1,\{\vartheta_{1,t}\}_{t\in\mathbb{R}}$ and $(\Omega,\mathcal{F},P,\{\vartheta_t\}_{t\in\mathbb{R}})$ are the requested
   parametric dynamical systems.

Given  $t\in\mathbb{R}, \omega=(\omega_1,\omega_2)\in\Omega$, we often write $\theta_t\omega=(\theta_t\omega_1,\theta_t\omega_2)$. Put
$$ z_1(\vartheta_t\omega_1)= -\lambda\int\limits_{-\infty}^0e^{\lambda s}(\vartheta_t\omega_1)(s)ds;\ \ \  z_2(\vartheta_t\omega_2)= -\sigma\int\limits_{-\infty}^0e^{\sigma s}(\vartheta_t\omega_2)(s)ds.$$
Then we have
$$dz_1(\vartheta_t\omega_1)+\lambda z_1(\vartheta_t\omega_1)dt=d\omega_1(t), \ \ dz_2(\vartheta_t\omega_1)+\sigma z_2(\vartheta_t\omega_1)dt=d\omega_2(t) .$$
By \cite{Arn} we know that both $ z_1(\vartheta_t\omega_1)$ and $ z_2(\vartheta_t\omega_2)$ are  continuous in $t$ on $\mathbb{R}$,  and further both $z_1(\omega_1)$ and $z_2(\omega_1)$
are  tempered.

Let $(\tilde{u},\tilde{v})$ satisfy  problem (1.1) and  write
\begin{align} \label{trans}
u(t,\tau,\omega,u_0)=\tilde{u}(t,\tau,\omega,\tilde{u}_0)-h_1z_1(\vartheta_t\omega_1)\ \  \mbox{and}\ v(t,\tau,\omega,v_0)=\tilde{v}(t,\tau,\omega,\tilde{v}_0)-h_2z_2(\vartheta_t\omega_2).
\end{align}
 Then $(u,v)$ solves the follow non-autonomous FitzHugh-Nagumo system
\begin{align} \label{pr1}
\frac{du}{dt}+\lambda u-\Delta u+\alpha v=f(x,u+h_1z_1(\vartheta_t\omega_1)+g(t,x)+\Delta z_1(\vartheta_t\omega_1)-\alpha h_2z_2(\vartheta_t\omega_2),
\end{align}
\begin{align} \label{pr2}
\frac{dv}{dt}+\sigma v-\beta u=h(t,x)+\beta h_1z_1(\vartheta_t\omega_1),
\end{align}
with the initial conditions  $u_0=\tilde{u}_0-h_1z_1(\vartheta_\tau\omega_1)$ and $v_0=\tilde{v}_0-h_2z_2(\vartheta_\tau\omega_2)$.

It is known (see \cite{Adili2}) that for every $(u_0,v_0)\in L^2(\mathbb{R}^N)\times  L^2(\mathbb{R}^N)$
 the problem (\ref{pr1})-(\ref{pr2})  possesses  a unique solution $(u,v)$ such that
$u\in C([\tau,+\infty),L^2(\mathbb{R}^N))\cap L^2(\tau,T, H^1(\mathbb{R}^N))\cap L^p(\tau,T, L^p(\mathbb{R}^N))$ and $v\in C([\tau,+\infty), L^2(\mathbb{R}^N))$. In addition,
the solution $(u,v)$ is continuous in $L^2(\mathbb{R}^N)\times L^2(\mathbb{R}^N)$ with respect to the initial value $(u_0,v_0)$ . Then formally $(\tilde{u},\tilde{v})=(u+h_1z_1(\vartheta_t\omega_1),v+h_2z_2(\vartheta_t\omega_2))$ is the solution to
problem (1.1) with the initial value  $\tilde{u}_0={u}_0+h_1z_1(\vartheta_\tau\omega_1)$ and $v\tilde{}_0={v}_0+h_2z_2(\vartheta_\tau\omega_2)$.

The continuous cocycle $\varphi$ over $\mathbb{R}$ and $(\Omega,\mathcal{F},P,\{\vartheta_t\}_{t\in\mathbb{R}})$  is defined by
\begin{align}\label{eq0}
& \ \ \ \ \ \ \ \ \ \varphi (t,\tau,\omega, (\tilde{u}_0,\tilde{v}_0))=(\tilde{u}(t+\tau,\tau,\vartheta_{-\tau}\omega,\tilde{u}_0),\tilde{v}(t+\tau,\tau,\vartheta_{-\tau}\omega,\tilde{v}_0))\notag\\
 &=(u(t+\tau,\tau,\vartheta_{-\tau}\omega,  u_0)+h_1z_1(\vartheta_t\omega_1),v(t+\tau,\tau,\vartheta_{-\tau}\omega,  v_0)+h_2z_2(\vartheta_t\omega_2) ),
\end{align} where $u_0=\tilde{u}_0-h_1z_1(\omega_1)$ and $v_0=\tilde{v}_0-h_2z_2(\omega_2)$.

Suppose that for every $\tau\in\mathbb{R}$ and $\omega\in \Omega$,
\begin{align} \label{D}
\lim\limits_{t\rightarrow+\infty}e^{-\delta t}\|D(\tau-t,\vartheta_{-t}\omega)\|_{L^2(\mathbb{R}^N)\times L^2(\mathbb{R}^N)}^2=0,
\end{align} where $\delta$ is as in (\ref{a5}).  Denote by $\mathcal{D}_\delta$ the collection of all families of bounded nonempty subsets of $L^2(\mathbb{R}^N)\times L^2(\mathbb{R}^N)$ which satisfies  (\ref{D}). Then it is obvious that $\mathcal{D}_\delta$ is inclusion closed.

In fact,
the existence  of
$\mathcal{D}_\delta$-pullback  attractors for the cocycle $\varphi$  defined by (\ref{eq0}) in the initial space $L^2(\mathbb{R}^N)\times L^2(\mathbb{R}^N)$  was proved in \cite{Adili2}.\\

\textbf{ Theorem 3.1(\cite{Adili2})}.\emph{ Assume that (\ref{a1})-(\ref{a5}) hold. Then the cocycle $\varphi$ has a unique $\mathcal{D}_\delta$-pullback attractor $\mathcal{A}=\{\mathcal{A}(\tau,\omega),\tau\in\mathbb{R}, \omega\in \Omega\}$ in $L^2(\mathbb{R}^N)\times L^2(\mathbb{R}^N)$. Here for every fixed $\tau\in\mathbb{R}$ and $\omega\in \Omega$, $\mathcal{A}(\tau,\omega)$ is given by
\begin{align}\label{L2}
\mathcal{A}(\tau,\omega)=\bigcap_{s>0}\overline{\bigcup_{t\geq s} \varphi(t,\tau-t,\vartheta_{-t}\omega, K(\tau-t,\vartheta_{-t}\omega))}^{L^2(\mathbb{R}^N)\times L^2(\mathbb{R}^N)},
\end{align} where $K$ is a closed and  measurable $\mathcal{D}_\delta$-pullback absorbing set of $\varphi$ in $L^2(\mathbb{R}^N)\times L^2(\mathbb{R}^N)$ which is stated by
\begin{align}\label{L2a}
 K(\tau,\omega)=\{(\tilde{u},\tilde{v})\in L^2(\mathbb{R}^N)\times L^2(\mathbb{R}^N);\|\tilde{u}\|^2+\|\tilde{v}\|^2\leq \rho(\tau,\omega)\},
\end{align} where $\rho(\tau,\omega)$ is the random constant:
\begin{align}\label{L2r}
 &\rho(\tau,\omega)=c(1+|z_1(\omega_1)|^2+|z_2(\omega_2)|^2)+ce^{-\delta \tau}\int\limits_{-\infty}^\tau e^{\delta s}(\|g(s,.)\|^2+\|h(s,.)\|^2)ds\notag\\
 &\ \  \ \ \ \ \  \ \ \ +c\int\limits_{-\infty}^0 e^{\delta s}(|z_1(\vartheta_{s}\omega_1)|^p+|z_2(\vartheta_{s}\omega_2)|^2)ds.
\end{align}Here $c$ is a deterministic positive constant.}

In the sequel, we will prove that  the set $\mathcal{A}$ defined by (\ref{L2}) is also a $\mathcal{D}_\delta$-pullback attractor in $L^\varpi(\mathbb{R}^N)\times L^2(\mathbb{R}^N)$ for
every $\varpi\in(2,p]$ with $p>2$.

\section{Existence of pullback random attractors in $L^{\varpi}(\mathbb{R}^N)\times L^{2}(\mathbb{R}^N)$}

 In this section, we will establish  the existence of $\mathcal{D}_\delta$-pullback attractors in $L^{\varpi}(\mathbb{R}^N)\times L^2(\mathbb{R}^N)(\varpi\in(2,p])$ for the cocycle $\varphi$ defined by (\ref{eq0}). Actually, we will show that the family of obtained
  pullback  attractors are the same one.

  Note that we do not increase any restrictions
 on the nonlinearity $f(x,s)$ and non-autonomous terms  $g$ and $h$ given in \cite{Adili2},
except that $\psi_1\in L^{p/2}(\mathbb{R}^N)$ as in (\ref{a1}).
For our problem,  the key point is prove that
the unbounded part of the first component $\tilde{u}$ of solutions of (1.1) vanishes in $L^p$, where $p>2$.

Hereafter, we denote by $\|.\|$ and $ \|.\|_p$ the norms in $L^2(\mathbb{R}^N)$ and $L^p(\mathbb{R}^N) (p>2)$, respectively.
Throughout this paper, the number $c$  is a generic positive constant independent of $\tau,\omega,D,\eta$ in any place.

Note that \cite{Adili2} has proved the following energy inequality (see the proof of Lemma 4.1 there):
\begin{align} \label{Borrow}
\frac{d}{dt}(\alpha\|v\|^2&+\beta\|u\|^2)+\delta(\alpha\|v\|^2+\beta\|u\|^2)+\frac{1}{2}\delta  \alpha\|v\|^2+\alpha_1\beta\|{u}\|_p^p\notag\\
&\leq \frac{4\beta}{\lambda}\|g(t,.)\|^2+\frac{4\alpha}{\sigma}\|h(t,.)\|^2+ c(|z_1(\vartheta_t\omega_1)|^{p}+|z_2(\vartheta_t\omega_2)|^2+1).
\end{align}

\subsection{Uniform estimates of solutions}

First we present a lemma which is sated in \cite{Adili2}.\\

\textbf{Lemma 4.1(\cite{Adili2}).}
 \emph{Assume that (\ref{a1})-(\ref{a5}) hold. Let $\tau\in\mathbb{R}, \omega\in\Omega$ and  $D=\{D(\tau,\omega);\tau\in\mathbb{R},\omega\in\Omega\}\in\mathcal{D}_\delta$. Then there exists a constant  $T=T(\tau,\omega,D)>0$
 such that for all $t\geq T$,  the solution $(u,v)$ of problem (\ref{pr1})-(\ref{pr2}) with $\omega$ replaced by $\vartheta_{-\tau}\omega$ satisfies
\begin{align} \label{Borrow1}
&\|u(\tau,\tau-t,\vartheta_{-\tau}\omega, \tilde{u}_0-h_1z_1(\vartheta_{-t}\omega_1)\|^2+\|v(\tau,\tau-t,\vartheta_{-\tau}\omega,\tilde{v}_0-h_2z_2(\vartheta_{-t}\omega_2)\|^2\leq R(\tau,\omega),
\end{align}
\begin{align} \label{Borrow2}
&\int\limits_{\tau-t}^\tau e^{\delta (s-\tau)}\Big(\|v(\tau,\tau-t,\vartheta_{-\tau}\omega, \tilde{v}_0-h_2z_2(\vartheta_{-t}\omega_2)\|^2+\|\tilde{u}(s,\tau-t,\vartheta_{-\tau}\omega, \tilde{u}_{0})\|_p^p\Big)ds\leq R(\tau,\omega),
\end{align}
where $(\tilde{u}_{0},\tilde{v}_{0})\in D_\delta(\tau-t,\vartheta_{-t}\omega)$ and
\begin{align} \label{radiu}
R(\tau,\omega)&=c+ce^{-\delta \tau}\int\limits_{-\infty}^\tau e^{\delta s}(\|g(s,.)\|^2+\|h(s,.)\|^2)ds\notag\\
&+c\int\limits_{-\infty}^{0} e^{\delta s }(|z_1(\vartheta_{s}\omega_1)|^{2p-2}+|z_1(\vartheta_{s}\omega_1)|^{p}+|z_1(\vartheta_{s}\omega_1)|^2+|z_2(\vartheta_{s}\omega_2)|^2)ds.
\end{align}}

In fact, for our problem we need further to show that (\ref{Borrow1}) holds in the  compact interval $[\tau-1.\tau]$. In particular, we have\\

\textbf{Lemma 4.2.}
 \emph{Assume that (\ref{a1})-(\ref{a5}) hold. Let $\tau\in\mathbb{R}, \omega\in\Omega$ and  $D=\{D(\tau,\omega);\tau\in\mathbb{R},\omega\in\Omega\}\in\mathcal{D}_\delta$. Then there exists a constant  $T=T(\tau,\omega,D)\geq2$
 such that for all $t\geq T$,  the first component of the solution of problem (\ref{pr1})-(\ref{pr2}) satisfies that for $\xi\in[\tau-1,\tau]$,
\begin{align} \label{}
\|u(\xi,\tau-t,\vartheta_{-\tau}\omega, \tilde{u}_0-h_1z_1(\vartheta_{-t}\omega_1)\|^2+\|{v}(\xi,\tau-t,\vartheta_{-\tau}\omega, \tilde{v}_{0}-h_2z_2(\vartheta_{-t}\omega_2)\|^2\leq R(\tau,\omega),\notag
\end{align}
where $(\tilde{u}_{0},\tilde{v}_{0})\in D_\delta(\tau-t,\vartheta_{-t}\omega)$ and $R(\tau,\omega)$ is as in (\ref{radiu})} only with a difference in the constant $c$.\\

\emph{Proof}\ \ We first let $t\geq2$ and $\xi\in[\tau-1,\tau]$.  Replacing $t$ by $s$ in (\ref{Borrow})
 and multiplying  by $e^{\delta(s-\tau)}$ and then integrating with respect to $s$ over $(\tau-t,\xi)$ we  deduce that,  along with $\omega$  replaced by $\vartheta_{-\tau}\omega$,
\begin{align} \label{Borrow11}
&e^{-\delta}(\beta\|u(\xi,\tau-t,\vartheta_{-\tau}\omega, \tilde{u}_0-h_1z_1(\vartheta_{-t}\omega_1)\|^2+\alpha\|v(\xi,\tau-t,\vartheta_{-\tau}\omega,\tilde{v}_0-h_2z_2(\vartheta_{-t}\omega_2)\|^2)\notag\\
&\ \ \ \ \ \ \  \ \ \ \ \  \leq \int\limits_{\tau-t}^\tau e^{\delta(s-\tau)}(\frac{4\beta}{\lambda}\|g(s,.)\|^2+\frac{4\alpha}{\sigma}\|h(s,.)\|^2)ds\notag\\
&\ \  \ \ \ \ \ \ \ \  \ \ \  \ + c\int\limits_{\tau-t}^\tau e^{\delta(s-\tau)} (|z_1(\vartheta_{s-\tau}\omega_1)|^{p}+|z_2(\vartheta_{s-\tau}\omega_2)|^2+1)ds+e^{-\delta t}(\beta\|u_0\|^2+\alpha\|v_0\|^2)\notag\\
&\ \ \ \ \ \ \  \ \ \ \ \ \  \leq ce^{-\delta\tau}\int\limits_{-\infty}^\tau e^{\delta s}(\|g(s,.)\|^2+\|h(s,.)\|^2)ds\notag\\
&\ \  \ \ \ \ \ \ \ \ \ \ \ \ \  + c\int\limits_{-\infty}^0 e^{\delta s} (|z_1(\vartheta_{s}\omega_1)|^{p}+|z_2(\vartheta_{s}\omega_2)|^2+1)ds
+ ce^{-\delta t}(\|\tilde{u}_0\|^2+\|\tilde{v}_0\|^2)\notag\\
&\ \ \ \ \ \ \  \ \ \ \ \ \ \ \ +ce^{-\delta t}(\|z_1(\vartheta_{-t}\omega_1)\|^2+\|z_2(\vartheta_{-t}\omega_2)\|^2).
\end{align} Note that  $z_1(\omega_1)$ and $z_2(\omega_2)$ are tempered and $(\tilde{u}_0,\tilde{v}_0)\in D(\tau-t,\vartheta_{-t}\omega)$. Then (\ref{Borrow11}) implies the desired.
 $\ \ \ \ \ \ \Box$\\

A result on the uniform $L^p$-estimate of the first component of  solution of problem (\ref{pr1})-(\ref{pr2}) on the compact interval $[\tau-1.\tau]$ is also needed  to show the uniform small of the unbounded part of solution in $L^p$ topology.\\

\textbf{Lemma 4.3.}
 \emph{Assume that (\ref{a1})-(\ref{a5}) hold. Let $\tau\in\mathbb{R}, \omega\in\Omega$ and  $D=\{D(\tau,\omega);\tau\in\mathbb{R},\omega\in\Omega\}\in\mathcal{D}_\delta$. Then there exists a constant  $T=T(\tau,\omega,D)\geq2$
 such that for all $t\geq T$,  the first component of the solution of problem (\ref{pr1})-(\ref{pr2}) satisfies that for $\xi\in[\tau-1,\tau]$,
\begin{align} \label{}
\|{u}(\xi,\tau-t,\vartheta_{-\tau}\omega, {u}_{0})\|_p^p\leq R(\tau,\omega),\notag
\end{align}
where $(\tilde{u}_{0},\tilde{v}_{0})\in D_\delta(\tau-t,\vartheta_{-t}\omega)$ and $R(\tau,\omega)$ is as in (\ref{radiu})} only with a difference in the constant $c$.\\

\emph{Proof}\ \
  Multiplying (\ref{pr1}) by $|u|^{p-2}u$ and then integrating over $\mathbb{R}^N$,  we have
\begin{align} \label{so10}
\frac{1}{p}\frac{d}{dt}\|u\|_p^p+\lambda\|u\|_p^p&\leq \int\limits_{\mathbb{R}^N}f(x,\tilde{u})|u|^{p-2}udx-\alpha \int\limits_{\mathbb{R}^N}v|u|^{p-1}udx\notag\\&+\int\limits_{\mathbb{R}^N}g(t,x)|u|^{p-2}udx+\int\limits_{\mathbb{R}^N}(\Delta h_1z_1(\vartheta_t\omega_1)-\alpha h_2z_2(\vartheta_t\omega_2))|u|^{p-2}udx.
\end{align}
By  using  (\ref{a1}) and (\ref{a2}) it infer us that
\begin{align} \label{}
&f(x,\tilde{u})u\leq -\frac{\alpha_1}{2^{p}}|{u}|^{p}+c|h_1z_1(\vartheta_t\omega_1)|^p+\psi_1(x)+\psi_2(x)|h_1z_1(\vartheta_t\omega_1)|,\notag
\end{align}
from which and by using Young inequality  we deduce that
\begin{align} \label{}
f(x,\tilde{u})u^{p-2}u
\leq-\frac{\alpha_1}{2^{p+1}}|{u}|^{2p-2}+\frac{\lambda}{2} u^p+c|h_1z_1(\vartheta_t\omega_1)|^{2p-2}+c\psi_1(x)^{p/2}+c|\psi_2(x)|^2.\notag
\end{align}
Then by $\psi\in L^{p/2}(\mathbb{R}^N$) and $\psi_2\in L^2(\mathbb{R}^N)$, we get the estimate of the nonlinearity  in (\ref{so10}),
\begin{align} \label{s011}
\int\limits_{\mathbb{R}^N} f(x,\tilde{u})u^{p-2}udx&\leq-\frac{\alpha_1}{2^{p+1}}\|u\|_{2p-2}^{2p-2}+\frac{\lambda}{2} \|u\|_p^p+c|z_1(\vartheta_t\omega_1)|^{2p-2}+c.
\end{align}
On the other hand, by using Young inequality again we get that
\begin{align} \label{so12}
&\alpha \int\limits_{\mathbb{R}^N}v|u|^{p-1}udx+\int\limits_{\mathbb{R}^N}g(t,x)|u|^{p-2}udx\notag\\
&+\int\limits_{\mathbb{R}^N}(\Delta h_1 z_1(\vartheta_t\omega_1)-\alpha h_2z_2(\vartheta_t\omega_2))|u|^{p-2}udx
\notag\\
&\leq \frac{\alpha_1}{2^{p+1}}\|{u}\|_{2p-2}^{2p-2}+c(\|v\|^2+\|g(t,.)\|^2+|z_1(\vartheta_t\omega_1)|^2+|z_2(\vartheta_t\omega_2)|^2).
\end{align}
Then combination (\ref{so10})-(\ref{so12}), we find that
\begin{align} \label{so14}
\frac{d}{dt}\|u\|_p^p+\delta\|u\|_p^p\leq c(\|v\|^2+\|g(t,.)\|^2+|z_1(\vartheta_t\omega_1)|^{2p-2}+|z_1(\vartheta_t\omega_1)|^2+|z_2(\vartheta_t\omega_2)|^2+1),
\end{align} where $\delta=\min\{\lambda,\sigma\}$. We first replace $t$ by $s$ in (\ref{so14}). Then by applying Gronwwall lemma (see also Lemma 5.1 in \cite{Zhao0}) over the interval $[\tau-t,\xi]$ for $t\geq 2$ and $\xi\in[\tau-1,\tau]$,
we find that, along with $\omega$ replaced by $\vartheta_{-\tau}\omega$,
\begin{align} \label{so15}
\|&u(\xi, \tau-t,\vartheta_{-\tau}\omega,{u}_0)\|_p^p\notag\\
&\leq \frac{e^{\delta}}{\xi-\tau+t}\int\limits_{\tau-t}^{\xi} e^{\delta (s-\tau)} \|u(s, \tau-t,\vartheta_{-\tau}\omega,\tilde{u}_0-h_1z_1(\vartheta_{-t}\omega_1)\|_p^pds\notag\\
&+ c\int\limits_{\tau-t}^{\xi} e^{\delta (s-\tau)}\|v(s,\tau-t,\vartheta_{-\tau}\omega,\tilde{v}_0-h_2z_2(\vartheta_{-t}\omega_2))\|^2ds+c\int\limits_{\tau-t}^{\xi} e^{\delta (s-\tau) } \|g(s,.)\|^2ds\notag\\
&+c\int\limits_{\tau-t}^{\xi} e^{\delta (s-\tau) }(|z_1(\vartheta_{s-\tau}\omega_1)|^{2p-2}+|z_1(\vartheta_{s-\tau}\omega_1)|^2+|z_2(\vartheta_{s-\tau}\omega_2)|^2+1)ds\notag\\
&\leq c\int\limits_{\tau-t}^{\tau} e^{\delta (s-\tau)} \|\tilde{u}(s, \tau-t,\vartheta_{-\tau}\omega,\tilde{u}_0)\|_p^pds\notag\\
&+c\int\limits_{\tau-t}^{\tau} e^{\delta (s-\tau)}\|v(s,\tau-t,\vartheta_{-\tau}\omega,\tilde{v}_0-h_2z_2(\vartheta_{-t}\omega_2))\|^2ds+ce^{-\delta \tau}\int\limits_{-\infty}^{\tau} e^{\delta s} \|g(s,.)\|^2ds\notag\\
&+c\int\limits_{-\infty}^{0} e^{\delta s }(|z_1(\vartheta_{s}\omega_1)|^{2p-2}+|z_1(\vartheta_{s}\omega_1)|^{p}+|z_1(\vartheta_{s}\omega_1)|^2+|z_2(\vartheta_{s}\omega_2)|^2)ds+c,
\end{align} where we have used that $\frac{1}{\xi-\tau+t}\leq 1$ for $t\geq2$. Then by applying (\ref{Borrow2}) to (\ref{so15}), we find that we have got the desired.
 $\ \ \ \ \ \ \Box$\\

\emph{Remark:} \ \  A similar calculation gives us that $u\in L^p(\mathbb{R}^N)$, which implies that for every $t>0,\tau\in\mathbb{R},\omega\in\Omega$ and $t>\tau$, the mapping $\varphi(t,\tau,\omega,.):L^2(\mathbb{R}^N)\times L^2(\mathbb{R}^N)\mapsto L^p(\mathbb{R}^N)\times L^2(\mathbb{R}^N)$. This demonstrate that the cocycle $\varphi$ defined by
(\ref{eq0}) satisfies condition \emph{(H1)}.\\

Let $u$ be the first component of solutions of problem (\ref{pr1})-(\ref{pr2}). Given
$\tau\in\mathbb{R},\omega\in\Omega$,
let $M=M(\tau,\omega)>0$ and
\begin{align} \label{}
\mathbb{R}^N(|u(\tau,\tau-t,\vartheta_{-\tau}\omega,u_0)|\geq M)=\{x\in \mathbb{R}^N; |u(\tau,\tau-t,\vartheta_{-\tau}\omega,u_0)|\geq M|\}.\notag
\end{align} Then we have\\

\textbf{Lemma 4.4.}
 \emph{Assume that (\ref{a1})-(\ref{a5}) hold. Let  $\tau\in\mathbb{R}, \omega\in\Omega$ and  $D=\{D(\tau,\omega);\tau\in\mathbb{R},\omega\in\Omega\}\in\mathcal{D}_\delta$. Then for any $\eta>0$,  there exist  constants $M=M(\tau,\omega,\eta)>0$ and $T=T(\tau,\omega,D)\geq2$
 such that for all $t\geq T$,  the first component $u$ of   solution $(u,v)$ of problem (\ref{pr1})-(\ref{pr2}) satisfies
\begin{align} \label{}
meas(\mathbb{R}^N(|u(\tau,\tau-t,\vartheta_{-\tau}\omega, \tilde{u}_0-h_1z_1(\omega_1)|\geq M)\leq\eta,\notag
\end{align} where meas(.) is the measure of a set and $(\tilde{u}_{0},\tilde{v}_{0})\in D(\tau-t,\vartheta_{-t}\omega)$.}\\

\emph{Proof}\ \  By (\ref{Borrow1}), there exist positive constants $T=T(\tau,\omega,D)>0$ and $R=R(\tau,\omega)$ such that for all $t\geq T$,
\begin{align} \label{Mes01}
\|u(\tau,\tau-t,\vartheta_{-\tau}\omega, \tilde{u}_0-h_1z_1(\vartheta_{-t}\omega_1)\|^2\leq R(\tau,\omega),
\end{align}where $(\tilde{u}_{0},\tilde{v}_{0})\in D(\tau-t,\vartheta_{-t}\omega)$.
Then for any $M>0$ and $|u(\tau)|=|u(\tau,\tau-t,\vartheta_{-\tau}\omega, \tilde{u}_0-h_1z_1(\omega_1)|\geq M$,  by  using (\ref{Mes01}) we have
\begin{align} \label{Mes02}
&M^2 meas(\mathbb{R}^N(|u(\tau,\tau-t,\vartheta_{-\tau}\omega, \tilde{u}_0-h_1z_1(\omega_1)|\geq M)\notag\\
&\leq \int\limits_{\mathbb{R}^N(|u(\tau)|\geq M)}|u(\tau,\tau-t,\vartheta_{-\tau}\omega, \tilde{u}_0-h_1z_1(\vartheta_{-t}\omega_1)|^2dx\leq R(\tau,\omega).
\end{align}By choosing $M=M(\tau,\omega,\eta)>(\frac{R(\tau,\omega)}{\eta})^{1/2}$, we get the result from (\ref{Mes02}). This ends the proof.
$\ \ \ \ \ \ \Box$\\

To prove the following lemma, we need to introduce some notations.
We will use a truncated version of the solution, $u-M$. Define
$$ (u-M)_+=\left\{
       \begin{array}{ll}
    u-M, \ \ \ \ \  \ \ \ \  \ u> M,\\
  0,\ \ \ \ \ \ \ \  \ \ \ \ \ \ \ \ \mbox{otherwise}.
       \end{array}
      \right.$$

The next lemma  will show that the unbounded part of the absolute value  $|u|$ vanishes  in $L^p$-norm on the state domain
$\mathbb{R}^N(|u(\tau,\tau-t,\vartheta_{-\tau}\omega,u_0)|\geq M)$ for $M$ large enough.\\

\textbf{Lemma 4.5.} Assume that (\ref{a1})-(\ref{a5}) hold. Let  $\tau\in\mathbb{R}, \omega\in\Omega$ and  $D=\{D(\tau,\omega);\tau\in\mathbb{R},\omega\in\Omega\}\in\mathcal{D}_\delta$.
 Then for any $\eta>0$, there exist  constants $M=M(\tau,\omega,\eta,D)>1$  and $T=T(\tau,\omega,D)\geq2$
such that the first component $\tilde{u}$ of solutions $(\tilde{u},\tilde{v})$ of problem (1.1) satisfies
\begin{align} \label{}
\sup_{t\geq T}\int\limits_{\mathbb{R}^N(|\tilde{u}(\tau)|\geq M)}|\tilde{u}(\tau,\tau-t,\vartheta_{-\tau}\omega,\tilde{u}_{0})|^p dx\leq \eta,\notag
\end{align} where $(\tilde{u}_{0},\tilde{v}_{0})\in D(\tau-t,\vartheta_{-t}\omega)$ and $\mathbb{R}^N(|\tilde{u}(\tau)|\geq M)=\mathbb{R}^N(|\tilde{u}(\tau,\tau-t,\vartheta_{-\tau}\omega,\tilde{u}_0)|\geq M)$.\\

\emph{Proof} \ \  Given $s\in[\tau-1,\tau]$, by replacing $\omega$ by $\vartheta_{-\tau}\omega$ in (\ref{pr1}) we get that
$$
u=u(s)=:u(s,\tau-t,\vartheta_{-\tau}\omega,u_0),v=v(s)=:v(s,\tau-t,\vartheta_{-\tau}\omega,v_0)\ \  s\in[\tau-1,\tau],
$$ satisfy the following  differential equations:
\begin{align} \label{p011}
\frac{du}{ds}+\lambda u-\Delta u+\alpha v=f(x,\tilde{u})+{g}(s,x)+\Delta h_1z_1(\vartheta_{s-\tau}\omega_1)-\alpha h_2z_2(\vartheta_{s-\tau}\omega_2),
\end{align}
where $(u_0,v_0)=(\tilde{u}_0-h_1z_1(\omega_1),\tilde{v}_0-h_2z_2(\omega_2))$ with $(\tilde{u}_0,\tilde{v}_0)\in D(\tau-t,\vartheta_{-t}\omega)$.

For fixed $\tau\in \mathbb{R}$ and $\omega\in\Omega$, we write
\begin{align} \label{p0M}
Z=Z(\tau,\omega)=\max\limits_{s\in[\tau-1,\tau]}\{\|h_1z_1(\vartheta_{s-\tau}\omega_1)\|_{L^\infty(\mathbb{R}^n)}\}, \ \ \  M=M(\tau,\omega)>\max\{1,Z\}.
 \end{align}
We  multiply (\ref{p011}) by  $(u-M)_+^{p-1}$  and integrate over $\mathbb{R}^N$ to get that for every $s\in[\tau-1,\tau]$,
\begin{align} \label{p01}
&\frac{1}{p}\frac{d}{ds}\int\limits_{\mathbb{R}^N} (u-M)_+^{p}dx+\lambda\int\limits_{\mathbb{R}^N} u(u-M)_+^{p-1}dx-\int\limits_{\mathbb{R}^N}\Delta u(u-M)_+^{p-1}dx\notag\\
&\ \  \ =-\alpha\int\limits_{\mathbb{R}^N}v(u-M)_+^{p-1}dx+\int\limits_{\mathbb{R}^N}f(x,\tilde{u})(u-M)_+^{p-1}dx\notag\\
&\ \ \  + \int\limits_{\mathbb{R}^N}{g}(s,x)(u-M)_+^{p-1}dx+\int\limits_{\mathbb{R}^N}(\Delta h_1z_1(\vartheta_{s-\tau}\omega_1)-\alpha h_2z_2(\vartheta_{s-\tau}\omega_2))(u-M)_+^{p-1}dx.
\end{align}
We now have to estimate every term in (\ref{p01}). First, it is obvious that
\begin{align} \label{p02}
-\int\limits_{\mathbb{R}^N}\Delta u(u-M)_+^{p-1}dx=(p-1)\int\limits_{\mathbb{R}^N}(u-M)_+^{p-2}|\nabla
u|^2dx\geq0,
\end{align}
\begin{align} \label{p03}
\lambda\int\limits_{\mathbb{R}^N} u(u-M)_+^{p-1}dx\geq \lambda\int\limits_{\mathbb{R}^N} (u-M)_+^{p}dx.
\end{align}
The most involved calculations are to estimate the nonlinearity in (\ref{p01}).
Note that
\begin{align} \label{p0p}
|\tilde{u}(s)|^{p-1}\geq 2^{2-p}|u(s)|^{p-1}-|h_1z_1(\vartheta_{s-\tau}\omega_1)|^{p-1},
\end{align}
and by (\ref{p0M}), if $u(s)>M$ then we have  for all $s\in [\tau-1,\tau]$,
\begin{align} \label{p034}
\tilde{u}(s)=u(s)+h_1z_1(\vartheta_{s-\tau}\omega_1)\geq u(s)-|h_1z_1(\vartheta_{s-\tau}\omega_1)|\geq u(s)-M>0.
\end{align}
Then by (\ref{a1}) and along with (\ref{p0p})-(\ref{p034}), we find that for $s\in[\tau-1,\tau]$ and  $u(s)>M$,
\begin{align} \label{p033}
f(x,\tilde{u})&\leq -\alpha_1|\tilde{u}|^{p-1}+\frac{1}{\tilde{u}}\psi_1(x)\leq -\alpha_12^{2-p}|u|^{p-1}
+\alpha_1|h_1z_1(\vartheta_{s-\tau}\omega_1)|^{p-1}+\frac{\psi_1(x)}{\tilde{u}}.
\end{align}
Then by (\ref{p033}) and along with (\ref{p034})
we get that for $u(s)>M$,
\begin{align} \label{}
 f(x,\tilde{u})&\leq
-\alpha_12^{2-p}|u|^{p-1}
+\alpha_1|h_1z_1(\vartheta_{s-\tau}\omega_1)|^{p-1}+|\psi_1(x)|(u-M)^{-1}\notag\\
&\leq -\alpha_12^{1-p}M^{p-2}(u-M)-\alpha_12^{1-p}(u-M)^{p-1}\notag\\
&+\alpha_1|h_1z_1(\vartheta_{s-\tau}\omega_1)|^{p-1}+|\psi_1(x)|(u-M)^{-1},\notag
\end{align}
from which it follows that
\begin{align} \label{p04}
&\ \ \   \int\limits_{\mathbb{R}^N}f(x,\tilde{u})(u-M)_+^{p-1}dx\notag\\&
\leq-\alpha_12^{1-p}M^{p-2}\int\limits_{\mathbb{R}^N}(u-M)_+^{p}dx-
\alpha_12^{1-p}\int\limits_{\mathbb{R}^N}(u-M)_+^{2p-2}dx\notag\\
&\ \  \  +\alpha_1\int\limits_{\mathbb{R}^N}|h_1z_1(\vartheta_{s-\tau}\omega_1)|^{p-1}(u-M)_+^{p-1}dx +\int\limits_{\mathbb{R}^N}|\psi_1(x)|(u-M)_+^{p-2}dx.
\end{align}
Since $h_1\in L^{2}(\mathbb{R}^N)\cap  L^{\infty}(\mathbb{R}^N)$, by the Sobolev interpolation, $h_1\in L^{2p-2}(\mathbb{R}^N)$. Then by using Young inequality
we have
\begin{align} \label{p041}
\alpha_1\int\limits_{\mathbb{R}^N}|h_1z_1(\vartheta_{s-\tau}\omega_1)|^{p-1}(u-M)_+^{p-1}dx\leq \frac{1}{2}\alpha_12^{1-p}\int\limits_{\mathbb{R}^N}(u-M)_+^{2p-2}dx+c|z_1(\vartheta_{s-\tau}\omega_1)|^{2p-2},
\end{align}
and by $\psi_1\in L^{p/2}(\mathbb{R}^N)$,
\begin{align} \label{p042}
\int\limits_{\mathbb{R}^N}|\psi_1(x)|(u-M)_+^{p-2}dx\leq \frac{1}{2}\lambda\int\limits_{\mathbb{R}^N}(u-M)_+^{p}dx+c\|\psi_1\|^{p/2}_{p/2}.
\end{align}
Combination (\ref{p04})-(\ref{p042}) we find
\begin{align} \label{p05}
&\ \ \ \ \   \int\limits_{\mathbb{R}^N}f(x,\tilde{u})(u-M)_+^{p-1}dx\notag\\&
\leq-\alpha_12^{1-p}M^{p-2}\int\limits_{\mathbb{R}^N}(u-M)_+^{p}dx-
\alpha_12^{-p}\int\limits_{\mathbb{R}^N}(u-M)_+^{2p-2}dx\notag\\
&\ \  \ \ \ +c|z_1(\vartheta_{s-\tau}\omega_1)|^{2p-2}+ \frac{1}{2}\lambda\int\limits_{\mathbb{R}^N}(u-M)_+^{p}dx+c\|\psi_1\|^{p/2}_{p/2}.
\end{align}
On the other hand, by using Young inequality, we also get that
\begin{align} \label{p051}
&-\alpha\int\limits_{\mathbb{R}^N}v(u-M)_+^{p-1}dx\leq \frac{1}{4}\alpha_12^{-p}\int\limits_{\mathbb{R}^N}(u-M)_+^{2p-2}dx
+c\int\limits_{\mathbb{R}^N(u(s)\geq M)}v^2dx,\\
&\int\limits_{\mathbb{R}^N}{g}(s,x)(u-M)_+^{p-1}dx\leq \frac{1}{8}\alpha_12^{-p}\int\limits_{\mathbb{R}^N}(u-M)_+^{2p-2}dx+c\int\limits_{\mathbb{R}^N(u\geq M)}{g}^2(s,x)dx,
\end{align}
and by $h_1\in H^2(\mathbb{R}^N)$ and $h_2\in L^2(\mathbb{R}^N)$,
\begin{align} \label{p052}
&\int\limits_{\mathbb{R}^N}(\Delta h_1z_1(\vartheta_{s-\tau}\omega_1)-\alpha h_2z_2(\vartheta_{s-\tau}\omega_2))(u-M)_+^{p-1}dx\notag\\
&\leq \frac{1}{8}\alpha_12^{-p}\int\limits_{\mathbb{R}^N}(u-M)_+^{2p-2}dx+c(|z_1(\vartheta_{s-\tau}\omega_1)|^2+|z_2(\vartheta_{s-\tau}\omega_2)|^2).
\end{align}
Then Combination (\ref{p01})-(\ref{p03}) and (\ref{p05})-(\ref{p052}), we obtain that
\begin{align} \label{p060}
\frac{d}{ds}&\int\limits_{\mathbb{R}^N} (u(s)-M)_+^{p}dx+
\alpha_12^{1-p}M^{p-2}\int\limits_{\mathbb{R}^N}(u(s)-M)_+^{p}dx\leq c(\|{g}(s,.)\|^2+\|v(s)\|^2)\notag\\
&\ \  \ \ +c(\|\psi_1\|_{p/2}^{p/2}+|z_1(\vartheta_{s-\tau}\omega_1)|^{2p-2}+|z_1(\vartheta_{s-\tau}\omega_1)|^2+|z_2(\vartheta_{s-\tau}\omega_2)|^2),
\end{align}
where the positive constant $c$  is independent of $\tau,\omega$ and $M$.
For convenience, we put
$$k=k(\tau,\omega,M)=2^{1-p}M^{p-2}, $$
and
$$\chi(\vartheta_{s}\omega)=c(\|\psi_1\|_{p/2}^{p/2}+|z_1(\vartheta_{s}\omega_1)|^{2p-2}+|z_1(\vartheta_{s}\omega_1)|^2+|z_2(\vartheta_{s}\omega_2)|^2).$$
  Then  we write (\ref{p060}) as
\begin{align} \label{p06}
\frac{d}{ds}\int\limits_{\mathbb{R}^N} (u(s)-M)_+^{p}dx&+
k\int\limits_{\mathbb{R}^N}(u(s)-M)_+^{p}dx\leq c(\|{g}(s,.)\|^2+\|v(s)\|^2)+\chi(\vartheta_{s-\tau}\omega),
\end{align} where $s\in[\tau-1,\tau]$.
Multiplying (\ref{p06}) by $e^{k(s-\tau)}$ and integrating first with respect to $s$
over $[\sigma,\tau]$  for $\sigma\leq \tau$  and then with respect to $\sigma$ over $[\tau-1,\tau]$, we
find that
 \begin{align} \label{p07}
&\int\limits_{\mathbb{R}^N} \Big(u(\tau,\tau-t,\vartheta_{-\tau}\omega,u_0)-M\Big)_+^{p}dx\notag\\
&\leq \int\limits_{\tau-1}^\tau e^{k(s-\tau)}\int\limits_{\mathbb{R}^N} \Big(u(s,\tau-t,\vartheta_{-\tau}\omega,u_0)-M\Big)_+^{p}dxds \notag\\
&+c\int\limits_{\tau-1}^\tau e^{k(s-\tau)}(\|g(s,.)\|^2+\|v(s,\tau-t,\vartheta_{-\tau}\omega,v_0)\|^2)ds +\int\limits_{\tau-1}^\tau e^{k(s-\tau)}\chi(\vartheta_{s-\tau}\omega)ds\notag\\
&\leq \int\limits_{\tau-1}^\tau e^{k(s-\tau)}(\|u(s,\tau-t,\vartheta_{-\tau}\omega,u_0)\|_p^{p}+\|v(s,\tau-t,\vartheta_{-\tau}\omega,v_0)\|^2)ds\notag\\
&+c\int\limits_{\tau-1}^\tau e^{k(s-\tau)}\|g(s,.)\|^2ds +\int\limits_{\tau-1}^\tau e^{k(s-\tau)}\chi(\vartheta_{s-\tau}\omega)ds\notag\\
&=I_1+I_2+I_3,
\end{align} where
$$
I_1=\int\limits_{\tau-1}^\tau e^{k(s-\tau)}(\|u(s,\tau-t,\vartheta_{-\tau}\omega,u_0)\|_p^{p}+\|v(s,\tau-t,\vartheta_{-\tau}\omega,v_0)\|^2)ds,
$$
$$
I_2=c\int\limits_{\tau-1}^\tau e^{k(s-\tau)}\|g(s,.)\|^2ds,
$$
and
$$
I_3=\int\limits_{\tau-1}^\tau e^{k(s-\tau)}\chi(\vartheta_{s-\tau}\omega)ds.
$$
We now calculate that $I_1,I_2$ and $I_3$ are  not larger than $\eta/3$ as $M\rightarrow+\infty$, respectively.
By Lemma 4.2 and 4.3, for every fixed $\tau\in\mathbb{R}$ and $\omega\in\Omega$  there exists $T_1=T_1(\tau,\omega,D)\geq 2$ such that for all $t\geq T_1$,
\begin{align} \label{p071}
&I_1\leq 2R(\tau,\omega)\int\limits_{\tau-1}^\tau e^{k(s-\tau)}ds\leq \frac{2R(\tau,\omega)}{k}
\leq\eta/3,
\end{align} as $M$ large enough, since $k\rightarrow +\infty$ as $M\rightarrow +\infty$. In order to estimate $I_2$, choosing $k>\delta$ and taking $\varsigma\in(0,1)$, we have
\begin{align} \label{}
I_2&=c\int\limits_{\tau-1}^{\tau-\varsigma} e^{k(s-\tau)}\|g(s,.)\|^2ds+c\int\limits_{\tau-\varsigma}^{\tau} e^{k(s-\tau)}\|g(s,.)\|^2ds\notag\\
&=ce^{-k\tau}\int\limits_{\tau-1}^{\tau-\varsigma} e^{(k-\delta)s}e^{\delta s}\|g(s,.)\|^2ds+ce^{-k\tau}\int\limits_{\tau-\varsigma}^{\tau} e^{ks}\|g(s,.)\|^2ds\notag\\
&\leq ce^{-k\varsigma} e^{\delta(\varsigma-\tau)}\int\limits_{-\infty}^{\tau} e^{\delta s}\|g(s,.)\|^2ds+c\int\limits_{\tau-\varsigma}^{\tau}\|g(s,.)\|^2ds.\notag
\end{align}
By (\ref{a5}), the first term above vanishes as $k\rightarrow+\infty$, and by $g\in L^2_{loc}(\mathbb{R}, L^2(\mathbb{R}^N))$  we can choose $\varsigma$ small enough such that the second term is small. Then when $M\rightarrow +\infty$, we have
\begin{align} \label{p07111}
&I_2\leq \eta/3.
\end{align}
To estimate $I_3$, we see that
\begin{align} \label{p07112}
&I_3=\int\limits_{-1}^0 e^{ks}\chi(\vartheta_{s}\omega)ds=\max_{0\leq s\leq1}\chi(\vartheta_{s}\omega)\frac{1}{k}\leq \eta/3,
\end{align} as $k$ large enough.
 Then from
(\ref{p07})-(\ref{p07112}), there exists  $M_1=M_1(\tau,\omega,\eta,D)>0$ large enough such that
\begin{align} \label{p10}
\sup_{t\geq T_1}\int\limits_{\mathbb{R}^N} \Big(u(\tau,\tau-t,\vartheta_{-\tau}\omega,u_0)-M_1\Big)_+^{p}dx&\leq \eta,
\end{align} where $T_1$ is as in (\ref{p071}).
If $u(\tau,\tau-t,\vartheta_{-\tau}\omega,u_0)\geq 2M_1$, then $u(\tau,\tau-t,\vartheta_{-\tau}\omega,u_0)-M_1\geq \frac{u(\tau,\tau-t,\vartheta_{-\tau}\omega,u_0)}{2},$ and therefore
by (\ref{p10}) it infer us  that
\begin{align} \label{p11}
&\sup_{t\geq T_1}\int\limits_{\mathbb{R}^N(u(\tau)\geq 2M_1)}
|{u}(\tau,\tau-t,\vartheta_{-\tau}\omega,{u}_0)|^{p}dx\notag\\
&\ \ \  \ \ \ \ \ \ \  \ \leq \sup_{t\geq T_1}2^{p}\int\limits_{\mathbb{R}^N} \Big(u(\tau,\tau-t,\vartheta_{-\tau}\omega,u_0)-M_1\Big)_+^{p}dx\leq  c\eta.
\end{align}
On the other hand, since by Sobolev interpolation we know that $h_1\in L^p(\mathbb{R}^N)$, then there is a small constant $\epsilon>0$ such that for any subset $E\subset \mathbb{R}^N$ with $meas (E)<\epsilon$,
\begin{align} \label{p1001}
\int\limits_E|h_1(x)|^pds\leq \frac{\eta}{|z_1(\omega_1)|^p}.
\end{align}
But by Lemma 4.4, there exist positive constants $M_2=M_2(\tau,\omega,\eta)>M_1$  and $T_2=T_2(\tau,\omega,D)>T_1$ such that
\begin{align} \label{p1002}
\sup_{t\geq T_2}meas(\mathbb{R}^N(|u(\tau,\tau-t,\vartheta_{-\tau}\omega,u_0)|\geq 2M_2))\leq \min\{\eta,\epsilon\}.
\end{align}
Then it follows from (\ref{p1001})-(\ref{p1002}) that
\begin{align} \label{p1003}
\sup_{t\geq T_2}\int\limits_{\mathbb{R}^N(u(\tau)\geq 2M_2)}|h_1z_1(\omega_1)|^pdx\leq \eta.
\end{align}
On the other hand, by (\ref{p11}) we get that
\begin{align} \label{p1101}
&\sup_{t\geq T_2}\int\limits_{\mathbb{R}^N(u(\tau)\geq 2M_2)}
|{u}(\tau,\tau-t,\vartheta_{-\tau}\omega,{u}_0)|^{p}dx\leq c\eta.
\end{align}
By taking advantage of the relation  $\tilde{u}(\tau,\tau-t,\vartheta_{-\tau}\omega,\tilde{u}_0)=u(\tau,\tau-t,\vartheta_{-\tau}\omega,u_0)+h_1z(\omega_1)$ and  association with (\ref{p0M}) we can deduce that
\begin{align} \label{}
\mathbb{R}^N(\tilde{u}(\tau,\tau-t,\vartheta_{-\tau}\omega,\tilde{u}_0)\geq 2M_2 +Z)\subseteq \mathbb{R}^N(u(\tau,\tau-t,\vartheta_{-\tau}\omega,u_0)\geq 2M_2),\notag
\end{align}
from which it implies that
\begin{align} \label{p80}
 &\sup_{t\geq T_2}\int\limits_{\mathbb{R}^N(\tilde{u}(\tau)\geq 2M_2+Z)}|\tilde{u}(\tau,\tau-t,\vartheta_{-\tau}\omega,\tilde{u}_0)|^{p}dx\notag\\
&\ \  \ \ \ \ \  \leq\sup_{t\geq T_2} 2^{p-1}\int\limits_{\mathbb{R}^N(u(\tau)\geq 2M_2)}|u(\tau,\tau-t,\vartheta_{-\tau}\omega,u_0)|^{p}dx\notag\\
&\ \ \ \ \ \ \ \ \ \ \ \ \  +\sup_{t\geq T_2}2^{p-1}\int\limits_{\mathbb{R}^N(u(\tau)\geq 2M_2)}|h_1z(\omega_1)|^pdx.
\end{align}
Therefore combination (\ref{p1003})-(\ref{p80}), we find that
\begin{align} \label{p801}
 \sup_{t\geq T_2}\int\limits_{\mathbb{R}^N(\tilde{u}(\tau)\geq 2M_2+Z)}
 |\tilde{u}(\tau,\tau-t,\vartheta_{-\tau}\omega,\tilde{u}_0)|^{p}dx\leq c\eta.
\end{align}
Similarly, we can deduce that there exists  $\tilde{M}_2=\tilde{M}_2(\tau,\omega,\eta,D)>0$ and $\tilde{T}_2=\tilde{T}_2(\tau,\omega,D)$ large enough such that
\begin{align} \label{p81}
\sup_{t\geq \tilde{T}_2}\int\limits_{\mathbb{R}^N(\tilde{u}(\tau)\leq - 2\tilde{M}_2-Z)}&|\tilde{u}(\tau,\tau-t,\vartheta_{-\tau}\omega,\tilde{u}_0)|^{p}dx\leq c\eta.
\end{align} Then (\ref{p801}) and (\ref{p81}) together imply the desired.
 $\ \ \  \ \ \ \ \ \ \Box$\\

\subsection{Existence of pullback random attractors in $L^{\varpi}(\mathbb{R}^N)\times L^2(\mathbb{R}^N)$}

 Note that the existence of pullback  attractor in $L^{2}(\mathbb{R}^N)\times L^2(\mathbb{R}^N)$  is prove in
\cite{Adili2}, which is stated by Theorem 3.1 in Section 3.
In this subsection, we will prove that this also happens  in $L^{\varpi}(\mathbb{R}^N)\times L^2(\mathbb{R}^N)$ for
every $\varpi\in(2,p]$. By Theorem 2.6, it remains to prove the asymptotic compactness of solutions of problem (1.1)
in $L^{\varpi}(\mathbb{R}^N)\times L^2(\mathbb{R}^N)$ for
every $\varpi\in(2,p]$, in which case we need to assume that $p>2$.

First the following lemma is borrowed from \cite{Adili2}.
\\

\textbf{Lemma 4.6.}  \emph{Assume that (\ref{a1})-(\ref{a5}) hold.   Then the  cocycle $\varphi$ associated with problem (1.1) is $\mathcal{D}_\delta$-asymptotically
compact in $L^2(\mathbb{R}^N)\times L^2(\mathbb{R}^N)$; that is,  for every  $\tau\in\mathbb{R}, \omega\in\Omega$,
the sequence $\{\varphi(t_n, \tau-t_n,\vartheta_{-t_n}\omega, (\tilde{u}_{0,n},\tilde{v}_{0,n}))\}$ has a convergent subsequence in
 $L^2(\mathbb{R}^N)\times L^2(\mathbb{R}^N)$ whenever $t_n\rightarrow+\infty$ and $(\tilde{u}_{0,n},\tilde{v}_{0,n})\in D(\tau-t_n,\vartheta_{-t_n}\omega)\in \mathcal{D}_\delta$.}
\\

By Proposition 2.8 and Lemma 4.5 and 4.6, we deduce the following, which reads\\

\textbf{Lemma 4.7.} \emph{Assume that (\ref{a1})-(\ref{a5}) hold. Then the cocycle $\varphi$ associated with
problem (1.1) is $\mathcal{D}_\delta$-asymptotically compact in $L^p(\mathbb{R}^N)\times L^2(\mathbb{R}^N)$; that is,  for every  $\tau\in\mathbb{R}, \omega\in\Omega$,
the sequence $\{\varphi(t_n, \tau-t_n,\vartheta_{-t_n}\omega, (\tilde{u}_{0,n},\tilde{v}_{0,n}))\}$ has a convergent subsequence in
 $L^p(\mathbb{R}^N)\times L^2(\mathbb{R}^N)$ whenever $t_n\rightarrow+\infty$ and $(\tilde{u}_{0,n},\tilde{v}_{0,n})\in D(\tau-t_n,\vartheta_{-t_n}\omega)\in \mathcal{D}_\delta$.}\\

Then by  Sobolev interpolation and Lemma 4.6 and Lemma 4.7,  we immediately  have\\

\textbf{Lemma 4.8.} \emph{Assume that (\ref{a1})-(\ref{a5}) hold. Then the  cocycle $\varphi$ associated with
problem (1.1) is $\mathcal{D}_\delta$-pullback asymptotically compact in $L^\varpi(\mathbb{R}^N)\times L^2(\mathbb{R}^N)$ for every $\varpi\in (2,p]$;
that is, for every  $\tau\in\mathbb{R}, \omega\in\Omega$, and $\varpi\in(2,p]$,
the sequence $\{\varphi(t_n, \tau-t_n,\vartheta_{-t_n}\omega, (\tilde{u}_{0,n},\tilde{v}_{0,n}))\}$ has a convergent subsequence in
 $L^\varpi(\mathbb{R}^N)\times L^2(\mathbb{R}^N)$ whenever $t_n\rightarrow+\infty$ and $(\tilde{u}_{0,n},\tilde{v}_{0,n})\in D(\tau-t_n,\vartheta_{-t_n}\omega)\in \mathcal{D}_\delta$.}\\

By Lemma 4.1, Lemma 4.6, Lemma 4.8, and  Theorem 2.6, we  get readily the following result, which implies that the obtained $\mathcal{D}_\delta$-pullback random attractor $\mathcal{A}$ in  $L^2(\mathbb{R}^N)\times L^2(\mathbb{R}^N)$ is compact and attracting  in  $L^\varpi(\mathbb{R}^N)\times L^2(\mathbb{R}^N)$ for every $\varpi\in(2,p]$ with $p>2$.\\

\textbf{Theorem 4.9.}  Assume that (\ref{a1})-(\ref{a5}) hold and $\mathcal{D}_\sigma$ is defined by (\ref{D}). Then the cocycle $\varphi$ associated with
problem (1.1)  possesses a unique $\mathcal{D}_\delta$-pullback random attractor $\mathcal{A}_\varpi=\{\mathcal{A}_\varpi(\tau,\omega);\tau\in\mathbb{R},\omega\in \Omega\}\in \mathcal{D}_\delta$ in $L^\varpi(\mathbb{R}^N)\times L^2(\mathbb{R}^N)$ for every $\varpi\in(2,p]$. Furthermore, for every $\tau\in\mathbb{R},\omega\in \Omega$,
\begin{align}\label{AA}
\mathcal{A}_\varpi(\tau,\omega)&=\cap_{s>0}\overline{\cup_{t\geq s}\varphi(t,\tau-t,\vartheta_{-t}\omega, K(\tau-t,\vartheta_{-t}\omega))}^{L^\varpi(\mathbb{R}^N)\times L^2(\mathbb{R}^N)}\\
&=\cap_{s>0}\overline{\cup_{t\geq s}\varphi(t,\tau-t,\vartheta_{-t}\omega, K(\tau-t,\vartheta_{-t}\omega))}^{L^2(\mathbb{R}^N)\times L^2(\mathbb{R}^N)}=\mathcal{A}(\tau,\omega),\notag
\end{align}
where $K=\{K(\tau,\omega);\tau\in\mathbb{R},\omega\in \Omega\}$ and $\mathcal{A}=\{\mathcal{A}(\tau,\omega);\tau\in\mathbb{R},\omega\in \Omega\}$ are the $\mathcal{D}_\delta$-pullback attractor and pullback absorbing set  of $\varphi$ in  $L^2(\mathbb{R}^N)\times L^2(\mathbb{R}^N)$  defined as in (\ref{L2}) and
(\ref{L2a}), respectively.\\

\textbf{Acknowledgments:}

 We are appreciated to Yangrong Li for his warm discussions and useful suggestions

  This work was supported by  Chongqing Basis and Frontier Research Project of China (No. cstc2014jcyjA00035) and National Natural Science Foundation of China
(No. 11271388).\\

\end{document}